\RequirePackage{fix-cm}
\documentclass[smallextended]{svjour3}       
\smartqed  
\usepackage{graphicx}
\usepackage{amssymb}
\usepackage{amsmath}
\usepackage{amsfonts}
\usepackage{graphicx}
\usepackage{graphics}
\usepackage{color}
\usepackage{latexsym}
\usepackage{url}

\input dmmacros.sty
\begin{document}

\title{Solving Heat Conduction Problems by
the Direct Meshless Local Petrov-Galerkin (DMLPG) method
}

\titlerunning{DMLPG for heat conduction}   

\author{Davoud Mirzaei \and Robert Schaback}


\institute{D. Mirzaei \at Department of Mathematics, University of Isfahan, 81745-163 Isfahan, Iran.\\
           \email{d.mirzaei@sci.ui.ac.ir}
           \and
           R. Schaback \at
           Institut f\"ur Numerische und Angewandte Mathematik, Universit\"at G\"ottingen, Lotzestra\ss{}e 16-18,
           D--37073 G\"ottingen, Germany.\\
           \email{schaback@math.uni-goettingen.de}
}

\date{Received: date / Accepted: date}

\maketitle

\begin{abstract}
As an improvement of the
{\em Meshless Local Petrov--Galerkin
(MLPG)}, the {\em Direct Meshless Local Petrov--Galerkin (DMLPG)} method
is applied here  to the
numerical solution of transient heat conduction problem.
The new technique is based on
{\em direct} recoveries of test functionals (local weak forms)
from values at nodes without any detour via
classical moving least squares (MLS) shape functions.
This leads to an absolutely cheaper scheme where the numerical integrations
will be done over low--degree polynomials rather than complicated MLS shape
functions. This eliminates the main disadvantage of MLS based methods in comparison with
finite element methods (FEM), namely the costs of numerical integration.


\end{abstract}

\section{Introduction}\label{SecIntro}
Meshless methods have received much attention in recent decades as
new tools
to overcome the difficulties of mesh generation and mesh refinement
in classical mesh-based methods such as the finite element method (FEM) and the
finite volume method (FVM).

The classification of numerical methods for solving PDEs should always start from
the classification of PDE problems themselves into {\em strong}, {\em weak},
or {\em local weak} forms. The first is the standard pointwise
formulation of differential equations and boundary conditions, the second is the
usual weak form dominating all FEM techniques, while the third
form splits the integrals of the usual global weak form into local
integrals over many small subdomains, performing the integration by parts on each
local integral. Local weak forms are the
basis of all variations of the {\em Meshless Local Petrov--Galerkin}
technique (MLPG) of S.N. Atluri and collaborators \cite{atluri:2005-1}. This
classification is dependent on the PDE problem itself,
and independent of numerical methods and the trial spaces used.
Note that these three formulations of the ``same'' PDE and boundary conditions
lead to three essentially different mathematical problems that cannot be
identified
and need a different mathematical analysis with respect to existence,
uniqueness, and stability of solutions.

Meshless {\em trial spaces}
mainly come via {\em Moving Least Squares} or {\em kernels} like
{\em Radial Basis Functions}.
They can consist of
{\em global} or {\em local} functions, but they should always
parametrize their trial functions ``{\it entirely in terms of nodes}''
\cite{belytschko-et-al:1996-1,wendland:2005-1} and require no triangulation or meshing.

A third classification
of PDE methods addresses where the discretization lives.
{\em Domain type} techniques work in the full global domain,
while {\em boundary type} methods work with exact solutions of the PDE and
just have to care for boundary conditions. This is independent of the other two
classifications.

Consequently, the literature should
confine the term ``meshless'' to be a feature of {\em trial spaces}, not of PDE problems
and their various formulations.
But many authors reserve the term {\em truly
  meshless} for meshless methods that
either do not require any discretization with a
background mesh for calculating
integrals or do not require integration at all.
These techniques have a great advantage in computational efficiency,
because numerical integration is the most time--consuming part in all
numerical methods based on local or global weak forms. This paper
focuses on a truly meshless method in this sense.

Most of the methods for solving PDEs in global weak form, such as the
Element-Free Galerkin (EFG) method \cite{belytschko-et-al:1994-1},
are not \emph{truly meshless} because
a triangulation is still required for
numerical integration.
The {\em Meshless Local Petrov-Galerkin} (MLPG) method
solves PDEs in local weak form and uses no global background mesh
to evaluate integrals because everything breaks down to some regular,
well-shaped and independent sub-domains.
Thus the MLPG is known as a truly meshless method.

We now focus on meshless methods using Moving Least Squares as trial functions.
If they solve PDEs in global or local weak form, they still suffer from the
cost of numerical integration. In these methods, numerical integrations are
traditionally done over
MLS shape functions and their derivatives.
Such shape functions are complicated and have no
closed form. To get accurate results, numerical
quadratures with many integration points are required. Thus the MLS subroutines
must be called very often, leading to high computational costs. In contrast to this,
the stiffness matrix in finite element methods (FEMs) is constructed by integrating
over polynomial basis functions which are much cheaper to evaluate. This relaxes
the cost of numerical integrations. For an account of the importance of
numerical integration within meshless methods, we refer the reader to
\cite{babuska-et-al:2009-1}.

To overcome this shortage within the MLPG based on MLS,
Mirzaei and Schaback \cite{mirzaei-schaback:2011-1} proposed a new technique,
\emph{Direct Meshless Local Petrov-Galerkin (DMLPG) method},
which avoids integration over MLS shape functions in MLPG and replaces
it by the much cheaper integration over polynomials. It ignores shape functions
completely. Altogether, the method is simpler, faster and often more accurate than the
original MLPG method. DMLPG uses a generalized MLS (GMLS)
method of \cite{mirzaei-et-al:2012-1}
which directly approximates boundary conditions and local weak forms as some
\emph{functionals}, shifting the numerical
integration into the MLS itself, rather than into an outside loop over calls to MLS
routines. Thus the concept of GMLS must be outlined first
in Section \ref{sect-gmls}
before we can go over to the DMLPG in Section \ref{sec-dmlpg}
and numerical results for heat conduction problems in Section \ref{sec-num}.

The analysis of heat conduction problems is important in engineering
and applied mathematics. Analytical solutions of heat equations
are restricted
to some special cases, simple geometries and specific boundary
conditions.
Hence, numerical methods are unavoidable.
Finite
element methods, finite volume methods, and
finite difference methods have been well
applied to transient heat analysis over the past few decades
\cite{minkowycz-et-al:1988-1}.
MLPG methods
were
also developed for heat transfer problems in many cases.
For instance,
J. Sladek et.al.
\cite{sladek-et-al:2003-1} proposed MLPG4 for transient heat
conduction analysis in functionally graded materials (FGMs) using Laplace transform techniques.
V. Sladek et.al.
\cite{sladek-et-al:2005-2} developed a local boundary
integral method for transient heat conduction in anisotropic and functionally graded media.
Both authors and their
collaborators
employed MLPG5 to analyze the heat
conduction in FGMs \cite{sladek-et-al:2004-1,sladek-et-al:2007-1}.

The aim of this paper is the
development
of DMLPG methods for heat conduction
problems.
This is the first time where DMLPG is applied to a time--dependent problem.
Moreover, compared to \cite{mirzaei-schaback:2011-1},
we will discuss all DMLPG methods, 
go into more details
and provide explicit formulae for the
numerical implementation.
DMLPG1/2/4/5 will be proposed, and the reason of ignoring DMLPG3/6
will be discussed. The new methods will be compared with the original MLPG
methods in a test problem, and then a problem in FGMs will be treated by DMLPG1.

In all application cases, the DMLPG method turned out to be superior to
the standard MLPG technique, and it provides excellent accuracy at low cost.
\section{Meshless methods and GMLS approximation}\label{sect-gmls}
Whatever the given PDE problem is and how it is discretized,
we have to find a function $u$ such that $M$ linear equations
\begin{equation}\label{eqlkubk}
\lambda_k(u)=\beta_k,\;1\leqslant k\leqslant M,
\end{equation}
defined by
$M$ linear {\em functionals} $\lambda_1,\ldots,\lambda_M$ and $M$
prescribed real values $\beta_1,\ldots,\beta_M$ are to be satisfied.
Note that weak formulations will involve functionals that integrate $u$ or a
derivative against some test function. The functionals can  discretize
either the  differential equation or some  boundary condition.

Now {\em meshless methods}
construct solutions from a {\em trial space} whose functions are
parametrized ``{\it entirely in terms of nodes}''
\cite{belytschko-et-al:1996-1}. We let these nodes
form a set $X:=\{x_1,\ldots,x_N\}$. Theoretically, meshless trial functions
can then be written as linear combinations of {\em shape functions}
$u_1,\ldots,u_N$ with or without the Lagrange conditions $u_j(x_k)=\delta_{jk},\;1\leqslant
j,k\leq N$ as
$$
  u(x)=\sum_{j=1}^Nu_j(x)  u(x_j)
$$
in terms of values at nodes, and this leads to solving the  system
(\ref{eqlkubk}) in the form
$$
\lambda_k(  u)
=\sum_{j=1}^N\lambda_k(u_j)  u(x_j)
=\beta_k,\;1\leqslant k\leqslant M
$$
approximately for the nodal values. Setting up the coefficient matrix
requires the evaluation of all functionals on all shape functions, and this
is a tedious procedure if the shape functions are not cheap to evaluate,
and it is even more tedious if the functionals consist of integrations
of derivatives against test functions.

But it is by no means mandatory to use shape functions at this stage at all.
If each functional $\lambda_k$ can be well approximated by a formula
\begin{equation}\label{eqlambdaapprox}
\lambda_k(  u)\approx \sum_{j=1}^N \alpha_{jk}   u(x_j)
\end{equation}
in terms of nodal values for smooth functions $  u$, the system
to be solved is
\begin{equation}\label{eqalphasys}
 \sum_{j=1}^N \alpha_{jk}   u(x_j)=\beta_k,\;1\leqslant k\leqslant M
\end{equation}
without any use of shape functions. There is no trial space, but
everything is still written in terms of values at nodes. Once the
approximate values $  u(x_j)$ at nodes are obtained, any
multivariate interpolation or approximation method can be used to generate
approximate values at other locations. This is a postprocessing step,
independent of PDE solving.

This calls for efficient ways to handle the approximations
(\ref{eqlambdaapprox})
to functionals in terms of nodal values.
We employ a generalized version of Moving Least Squares (MLS),
adapted from \cite{mirzaei-et-al:2012-1}, and without using shape functions.

The techniques of \cite{mirzaei-et-al:2012-1} and \cite{mirzaei-schaback:2011-1}
allow to calculate coefficients $\alpha_{jk}$ for (\ref{eqlambdaapprox})
very effectively as follows. We fix $k$ and consider just $\lambda:=\lambda_k$.
Furthermore, the set $X$  will be formally replaced
by a much smaller subset that consists only of the
nodes that are locally necessary to calculate a good approximation of
$\lambda_k$,
but we shall keep $X$ and $N$ in the notation. This reduction of the
node set for the approximation of $\lambda_k$ will ensure sparsity
of the final coefficient matrix in (\ref{eqalphasys}).

Now we have to calculate a coefficient vector
$a(\lambda_k)=(\alpha_{1k},\ldots,\alpha_{Nk})^T\in\R^N$
for (\ref{eqlambdaapprox}) in case of $\lambda=\lambda_k$. We choose  a
space $\calp$ of polynomials which is large enough to let
zero be the only polynomial $p$ in $\calp$ that vanishes on $X$.
Consequently, the
dimension $Q$ of $\calp$ satisfies $Q\leqslant N$,
and the $Q\times N$ matrix
 $P$ of values $p_i(x_j)$ of a basis $p_1,\ldots,p_Q$ of $\calp$
has rank $Q$.
Then for any  vector
$w=(w_1,\ldots,w_N)^T$ of positive weights,
the generalized MLS solution $a(\lambda)$ to (\ref{eqlambdaapprox}) can be written as
\begin{equation}\label{aWl}
a(\lambda_k)=WP^T(P\,W\,P^T)^{-1}\lambda_k(\calp)
\end{equation}
where $W$ is the diagonal matrix with diagonal $w$ and
$\lambda_k(\calp)\in\R^Q$ is the vector with values
$\lambda_k(p_1),\ldots,\lambda_k(p_Q)$.

Thus it suffices to evaluate $\lambda_k$ on low--order polynomials, and since
the coefficient matrix in (\ref{aWl}) is independent of $k$,
one can use the same matrix for different $\lambda_k$ as long as $X$ does not
change locally. This will significantly speed up numerical calculations, if the functional
$\lambda_k$ is complicated, e.g. a numerical integration against a test function.
Note that the MLS is just behind the scene, no shape functions occur. But the
weights will be defined locally in the same way as in the usual MLS, e.g.
we choose a continuous function $\phi:[0,\infty)\rightarrow[0,\infty)$
with
\begin{itemize}
\item $\phi(r)> 0, ~ 0\leqslant r< 1,$
\item $\phi(r)=0,~ r\geqslant 1,$
\end{itemize}
and define
\begin{equation*}\label{weightw}
w_j(x)=\phi\left( \frac{\|x-x_j\|_2}{\delta}\right)
\end{equation*}
for $\delta>0$ as a weight function, if we work locally near a point $x$.
\vspace{-6pt}
\section{MLPG Formulation of Heat Conduction}\label{sec-mlpg}
\vspace{-2pt}
In the Cartesian
coordinate system, the transient temperature field in a heterogeneous
isotropic medium is governed by the diffusion equation
\begin{equation}\label{govern1}
\rho(x) c(x)\frac{\partial u}{\partial t}(x,t) = \nabla \cdot (\kappa \nabla u) +f(x,t),
\end{equation}
where $x\in\Omega$ and $0\leqslant t\leqslant t_F$ denote the
space and time
variables, respectively, and $t_F$ is the
final time. The initial and boundary conditions are
\begin{align}
u(x,0)&=u_0(x),\quad x\in \Omega,\label{initialcond}\\
u(x,t)&=u_D(x,t),\quad x\in \Gamma_D,
\quad 0\leqslant t\leqslant t_F, \label{drichletcond}\\
\kappa(x)\frac{\partial u}{\partial n}(x,t)
&=u_N(x,t),\quad x\in \Gamma_N,\quad 0\leqslant t\leqslant t_F\label{neumanncond}.
\end{align}
In (\ref{govern1})-(\ref{neumanncond}),
$u(x,t)$ is the temperature field, $\kappa(x)$ is the thermal conductivity dependent on
the
spatial
variable $x$, $\rho(x)$ is the mass density and $c(x)$ is the
specific heat, and
$f(x,t)$ stands for the internal heat source
generated per unit volume. Moreover,
$n$ is the unit
outward normal to the boundary $\Gamma$, $u_D$
and $u_N$ are specified
values on the Dirichlet boundary $\Gamma_D$ and
Neumann boundary $\Gamma_N$ where $\Gamma=\Gamma_D\cup\Gamma_N$.

Meshless methods
write everything entirely in terms of scattered
nodes forming a set $X = \{x_1,x_2,\ldots, x_N\}$
located in the spatial domain $\Omega$ and its boundary $\Gamma$.
In the standard MLPG, around each $x_k$ a small subdomain $\Omega_s^k\subset
\overline{\Omega}=\Omega\cup \Gamma$ is chosen such that integrations over
$\Omega_s^k$ are comparatively cheap. For instance, $\Omega_s^k$
is conveniently taken to be the intersection of
$\overline \Omega$ with a ball $B(x_k,r_0)$ of radius $r_0$ or a cube
(or a square in 2D)
$S(x_k,r_0)$ centered at $x_k$ with side-length $r_0$.
On these subdomains, the PDE including boundary conditions is stated in
a localized weak form
\begin{equation}\label{eq-ulwf-v}
\frac{\partial}{\partial t}\int_{\Omega_s^{k}} \rho c u v\, d\Omega =
\int_{\Omega_s^{k}} \nabla \cdot (\kappa \nabla u) v \, d\Omega + \int_{\Omega_s^{k}}f v \, d\Omega,
\end{equation}
for an appropriate {\em test function} $v$.
Applying integration by parts, this weak equation can be partially
symmetrized to become
the {\em first} local weak form
\begin{equation}\label{eq-lwf-v}
\frac{\partial}{\partial t}\int_{\Omega_s^{k}} \rho c u v\, d\Omega
= \int_{\partial \Omega_s^{k}} \kappa \frac{\partial u}{\partial n}\, v \, d
\Gamma
-\int_{\Omega_s^{k}}\kappa \nabla u\cdot \nabla v \, d\Omega + \int_{\Omega_s^{k}}f v \, d\Omega.
\end{equation}
The {\em second}
local weak form, after rearrangement of \eqref{govern1} and integration
by parts twice, can be obtained as
\begin{equation}\label{eq-2ulwf-v}
\begin{split}
\frac{\partial}{\partial t}\int_{\Omega_s^{k}} \frac{1}{\kappa}\rho c u v \,
d\Omega &
= \int_{\Omega_s^k} u \Delta v\, d\Omega -
\int_{\partial \Omega_s^k} u \frac{\partial v}{\partial n}\, d\Gamma
+ \int_{\partial \Omega_s^k} v \frac{\partial u}{\partial n}\, d\Gamma \\
& + \int_{\Omega_s^{k}}\frac{1}{\kappa}\nabla\kappa\cdot \nabla u\, v  \,
d\Omega
+ \int_{\Omega_s^{k}}\frac{1}{\kappa}f v \, d\Omega.
\end{split}
\end{equation}
If the boundary of the local domain $\Omega_s^k$
hits the boundary of $\Omega$, the MLPG inserts boundary data at the appropriate
places in order to care for boundary conditions.
Since these local weak equations are all affine--linear in $u$ even after
insertion of boundary data, the equations of MLPG are all of the form
(\ref{eqlkubk})
after some rearrangement, employing certain linear functionals $\lambda_k$.
In all cases, the MLPG evaluates these functionals on shape functions, while
our DMLPG method will use the GMLS approximation of Section
\ref{sect-gmls}  without any shape function.

However, different choices of test functions $v$ lead to the six different
well--known types of MLPG.
The variants MLPG1/5/6 are based on the weak formulation \eqref{eq-lwf-v}.
If $v$ is chosen such that the first integral in the right hand side
of (\ref{eq-lwf-v})
vanishes, we have MLPG1.
In this case $v$ should vanish on $\partial \Omega_s^k$. If the Heaviside
step function $v$ on local domains is used as test function,
the second integral disappears and we have a pure
local boundary integral form in the right hand side.
This is MLPG5.
In MLPG6, the trial and test functions come from the same space.
MLPG2/3 are based on the local unsymmetric weak formulation
\eqref{eq-ulwf-v}. MLPG2 employs Dirac's delta
function as the test function in each  $\Omega_s^k$,
which leads to a pure collocation method. MLPG3 employs the error
function as the test function in each $\Omega_s^k$. In this method,
the test functions can be the same as for the discrete least squares method.
The test functions and the trial functions come from the same space in MLPG3.
Finally, MLPG4 (or LBIE) is based on the weak form \eqref{eq-2ulwf-v},
and a modified fundamental solution of the
corresponding elliptic spatial equation is employed as a test function
in each subdomain.

We describe these types in more detail later,
along with the way we modify them
when going from MLPG to DMLPG.
\section{DMLPG Formulations}\label{sec-dmlpg}
Independent of which variation of MLPG we go for,
the DMLPG has its special ways to handle boundary conditions,
and we describe these first.

Neither Lagrange multipliers
nor penalty parameters are introduced
into the local weak forms, because the Dirichlet boundary conditions are imposed
directly.
For nodes $x_k\in \Gamma_D$, the values
$u(x_k,t)=u_D(x_k,t)$ are known from the Dirichlet boundary conditions.
To connect them properly to nodal values $u(x_j,t)$
in neighboring points $x_j$ inside the domain
or on the Neumann boundary, we turn the GMLS philosophy upside down
and ask for coefficients $a_j(x_k)$ that allow to reconstruct nodal values at
$x_k$ from nodal values at the $x_j$. This amounts to setting
$\lambda_k=\delta_{x_k}$ in Section \ref{sect-gmls} , and we get
localized equations for Dirichlet boundary points $x_k$ as
\begin{equation}\label{dirichlet-impose}
\sum_{j=1}^N a_j(x_{k}) u(x_j,t) = u_D (x_{k},t),
\quad x_{k}\in\Gamma_D, \quad t\in [0,t_F].
\end{equation}
Note that the coefficients are time--independent.
In matrix form,
(\ref{dirichlet-impose})
can be written as
\begin{equation}\label{dirichlet-sys}
B \u(t) = \u_D(t),
\end{equation}
where
$
\u(t)\in\R^N
$
is
the
time--dependent vector of nodal values at $x_1,x_2,...,x_N$.
These equations are added into
the full matrix setup at the appropriate places, and they are in truly meshless
form, since they involve only values at nodes and are without numerical
integration. Note that (\ref{eq-lwf-v}) has no integrals
over the Dirichlet boundary, and thus we can impose Dirichlet conditions
always in the above strong form.
For \eqref{eq-2ulwf-v} there are two possibilities.
We can impose the Dirichlet boundary conditions either in the local weak form or
in the collocation form \eqref{dirichlet-impose}. Of course the
latter is the cheaper one.

We now turn to Neumann boundary conditions.
They can be imposed in the same way as
Dirichlet boundary conditions by assuming
$\lambda_k(u) = \frac{\partial u}{\partial n}(x_k)$
in the
GMLS approximation
\begin{equation}\label{neumann-impose}
\sum_{j=1}^N a_j(x_{k}) u(x_j,t) = \frac{\partial u}{\partial n}(x_{k},t),
\quad x_{k}\in\Gamma_N, \quad t\in [0,t_F].
\end{equation}
Note that the coefficients again are time--independent,
and we get a linear system like (\ref{dirichlet-sys}),
but with a vector $\u_N(t)$ of nodal values of
normal derivatives in the right--hand side. This is collocation as in
subsection \ref{subsec-dmlpg2}. But it is often more accurate to impose
Neumann conditions directly into the local weak forms
\eqref{eq-lwf-v} and \eqref{eq-2ulwf-v}. We will describe this
in more detail in the following subsections.
We now turn the different variations of the MLPG method into variations of the DLMPG.


\subsection{DMLPG1/5}
These methods are based on the local weak form \eqref{eq-lwf-v}. This form recasts to
\begin{equation}\label{eq-lwf1}
\begin{split}
\frac{\partial}{\partial t}\int_{\Omega_s^{k}} \rho c u v\, d\Omega
&+\int_{\Omega_s^{k}}\kappa \nabla u\cdot \nabla v \, d\Omega
-\int_{\partial\Omega_s^k\setminus \Gamma_N}
\kappa \frac{\partial u}{\partial n}v\, d\Gamma\\
& =\int_{\Gamma_N\cap \partial\Omega_s^k} u_N v\, d\Gamma
+ \int_{\Omega_s^{k}}f v \, d\Omega
\end{split}
\end{equation}
after inserting the Neumann boundary data
from {\eqref{neumanncond}}, when
the domain $\Omega_s^k$ of (\ref{eq-lwf-v})
hits the Neumann boundary $\Gamma_N$.
All integrals in the top part
of \eqref{eq-lwf1}
can be efficiently approximated
by GMLS approximation
of Section \ref{sect-gmls}
as purely spatial formulas
\begin{equation}\label{eqlamjk}
\begin{array}{rclcl}
\lambda_{1,k}(u)
&:=&
\displaystyle{\int_{\Omega_s^{k}} \rho c u v\, d\Omega}
&\approx&\displaystyle{
  \widehat {\lambda_{1,k}(u)}= \sum_{j=1}^N a_{1,j} (x_k)u(x_j),}\\
\lambda_{2,k}(u)
&:=&\displaystyle{  -\int_{\Omega_s^{k}}\kappa \nabla u\cdot \nabla v \, d\Omega}
&\approx&\displaystyle{
\widehat {\lambda_{2,k}(u)}= \sum_{j=1}^N a_{2,j} (x_k)u(x_j),}\\
\lambda_{3,k}(u)
&:=& \displaystyle{ -\int_{\partial\Omega_s^k\setminus \Gamma_N}
\kappa \frac{\partial u}{\partial n}v\, d\Gamma}
&\approx&\displaystyle{
\widehat {\lambda_{3,k}(u)}= \sum_{j=1}^N a_{3,j} (x_k)u(x_j)}.
\end{array}
\end{equation}

While the two others can always be summed up,
the first formula, if applied to time--varying functions, has to be modified
into
$$
\displaystyle{\frac{\partial}{\partial t}\int_{\Omega_s^{k}} \rho c u v\, d\Omega}
\approx\displaystyle{
 \sum_{j=1}^N a_{1,j} (x_k)\frac{\partial}{\partial t}u(x_j,t)}
$$
and expresses the main PDE term not in terms of values at nodes, but rather in
terms of time derivatives of values at nodes.

Again, everything is expressed in terms of values at nodes,
and the coefficients are time--independent. Furthermore,
Section \ref{sect-gmls}  shows that the $u$ part of the
integration runs over low--order polynomials, not over any shape functions.

The third functional can be omitted if the
test function $v$ vanishes on $\partial\Omega_s^k\setminus \Gamma_N$. This is DMLPG1.
An example of such a
test function is
\begin{equation*}\label{weight}
v=v(x;x_k) = \phi\left( \frac{\|x-x_k\|_2}{r_0}\right),
\end{equation*}
where $\phi$ is the weight function in the MLS approximation
with the radius $\delta$ of the
support of the weight function being replaced by the radius
$r_0$ of the local domain $\Omega_s^k$.

In DMLPG5, the local test function is the constant $v=1$.
Thus the functionals $\lambda_{2,k}$ of (\ref{eqlamjk}) are not needed, and the
integrals for $\lambda_{1,k}$ take a simple form, if $c$ and $\rho$ are simple.
DMLPG5 is slightly cheaper than DMLPG1,
because the domain integrals of $\lambda_{2,k}$  are replaced by the boundary
integrals of $\lambda_{3,k}$.

Depending on
which parts of the functionals are present or not, we finally get
a time--dependent system of the form
\begin{equation}\label{eqfullsys}
A^{(1)} \frac{\partial}{\partial t} \u(t) +A^{(\ell)} \u(t) =\b(t),
\quad \ell=2\,\, \mbox{or}\,\, 3
\end{equation}
where $\u(t)$  is the time--dependent vector
$$
\u(t)=(u(x_1,t),\ldots,u(x_N,t))^T\in\R^N
$$
of nodal values,
$\b(t)\in\R^M$ collects the time--dependent right--hand sides
with components
$$
b_k = \int_{\Omega_s^{k}}f(x,t)v(x;x_k)  \, d\Omega +\int_{\Gamma_N\cap \partial\Omega_s^k}u_N(x,t) v(x;x_k)\, d\Gamma,
$$
and $A_{kj}^{(\ell)} = a_{\ell,j}(x_k) $, $\ell=1,2,3$. The $k$-th row of $A^{(\ell)}$ is
\begin{align*}
\a_k^{(\mathrm{\ell})}  &= WP (P WP^T)^{-1} \lambda_{\ell,k}(\calp),\quad \ell = 1,2,3,
\end{align*}
where
\begin{align*}
\lambda_{1,k}(\calp)=& \left[\int_{\Omega_s^{k}} \rho c p_1 v\, d\Omega,
\int_{\Omega_s^{k}} \rho c p_2 v \, d\Omega, \ldots ,
\int_{\Omega_s^{k}} \rho c p_Q v \, d\Omega \right]^T, \\ 
\lambda_{2,k}(\calp)=& -\left[\int_{\Omega_s^{k}}\kappa \nabla p_1\cdot \nabla v
  \, d\Omega,
\int_{\Omega_s^{k}}\kappa \nabla p_2\cdot \nabla v \, d\Omega, \ldots,
\int_{\Omega_s^{k}}\kappa \nabla p_Q\cdot \nabla v \, d\Omega\right]^T,\\ 
\lambda_{3,k}(\calp)=& \left[
\int_{\partial\Omega_s^k\setminus \Gamma_N} \kappa \frac{\partial p_1}{\partial
  n}v\, d \Gamma,
\int_{\partial\Omega_s^k\setminus \Gamma_N} \kappa \frac{\partial p_2}{\partial
  n}v\, d \Gamma, \ldots,
\int_{\partial\Omega_s^k\setminus \Gamma_N} \kappa \frac{\partial p_Q}{\partial n}v  \, d \Gamma\right]^T.
\end{align*}
As we can immediately see, {\em numerical integrations are done over
low-degree polynomials $p_1, p_2, ..., p_Q$
only, and no shape function is needed at all. This reduces
the cost of numerical integration in MLPG methods significantly.}


\subsection{DMLPG2}\label{subsec-dmlpg2}
In this method, the test function $v$ on the local domain
$\Omega_s^k$ in \eqref{eq-ulwf-v} is replaced by the test functional $\delta_{x_k}$, i.e.
we have strong collocation of the PDE and all boundary conditions.
Depending on where $x_k$ lies, one can have the functionals
\begin{equation}\label{eqallmu}
\begin{array}{rcl}
\mu_{1,k}(u)&:=& u(x_k),\\
\mu_{2,k}(u)&:=& \displaystyle{ \frac{\partial u}{\partial n}(\rho c u)(x_k)},\\
\mu_{3,k}(u)&:=& \nabla\cdot (\kappa \nabla u)(x_k)\\
\end{array}
\end{equation}
connecting $u$ to Dirichlet, Neumann, or PDE data. The first form is used
on the Dirichlet boundary, and leads to   (\ref{dirichlet-impose})
and (\ref{dirichlet-sys}). The second applies to points on the Neumann boundary
and is handled by (\ref{neumann-impose}),
while the third can occur anywhere in $\overline \Omega$ independent of the other possibilities.
In all cases, the GMLS method of Section \ref{sect-gmls} leads to
approximations of the form
$$
\mu_{i,k}(u)\approx \sum_{j=1}^N a_{i,j}(x_k)u(x_j),\;i=1,2,3
$$
entirely in terms of nodes, where values on nodes on the Dirichlet boundary can
be replaced by given data.

This DMLPG2 technique is a pure  collocation method and requires
no numerical integration at all.
Hence it is
truly meshless and
the
cheapest
among all versions of DMLPG and MLPG.
But
it needs
higher order derivatives,
and thus the order of convergence is reduced
by the order of the derivative taken.
Sometimes DMLPG2 is called {\em Direct MLS Collocation (DMLSC) method} \cite{mirzaei-schaback:2011-1}.

It is worthy to note that the recovery of a
functional such as $\mu_{2,k}(u)$ or $\mu_{3,k}(u)$  in (\ref{eqallmu})
using GMLS approximation gives {\em GMLS derivative approximation}.
These kinds of derivatives have
been comprehensively investigated in
\cite{mirzaei-et-al:2012-1} and a rigorous error bound
was derived for them. Sometimes they are called {\em diffuse} or
{\em uncertain} derivatives, because they are not derivatives of shape functions,
but \cite{mirzaei-et-al:2012-1} proves there is
nothing diffuse or uncertain about them and they are
direct and usually very good numerical approximation of corresponding function derivatives.


\subsection{DMLPG4}
This method is based on the local weak form \eqref{eq-2ulwf-v} and uses the
{\em fundamental solution} of the corresponding elliptic spatial equation as test
function. Here we describe it for a two--dimensional problem.
To reduce the unknown quantities in local weak
forms, the concept of {\em companion solutions}
was introduced in \cite{zhu-et-al:1998-1}.
The companion solution of a 2D Laplace operator is
\begin{equation*}
v(x;y) = \frac{1}{2\pi} \ln \frac{r_0}{r}, \quad r= \|x-y\|_2,
\end{equation*}
which corresponds to
the Poisson equation $\Delta v(x;y) +\delta(r) = 0$
and thus is a {\em fundamental} solution vanishing for $r=r_0$.
Dirichlet boundary conditions for DMLPG4
are imposed as in \eqref{dirichlet-impose}.
The resulting local integral equation corresponding to
a node $x_k$ located inside the domain
or on the Neumann part of the boundary is
\begin{equation}\label{eq-ulwf2-v}
\begin{split}
\frac{\partial}{\partial t}\int_{\Omega_s^{k}} \frac{1}{\kappa}\rho c u v \, d\Omega
- & \alpha_k u(x_k) +  \dashint_{\partial \Omega_s^{k}} \frac{\partial v}{\partial n}u  \, d \Gamma -
\int_{\Omega_s^{k}}\frac{1}{\kappa}\nabla\kappa\cdot \nabla u\, v  \, d\Omega \\
&= \int_{\partial \Omega_s^k\cap\Gamma_N}u_N v \, d\Gamma
+ \int_{\Omega_s^{k}}\frac{1}{\kappa}f v \, d\Omega,
\end{split}
\end{equation}
where $\alpha_k$
is a coefficient that depends
on where the source
point $x_k$ lies.
It is $1/2$ on the smooth
boundary,
and $\theta_k/(2\pi)$
at a corner where the interior angle at the point $x_k$ is $\theta_k$. The symbol
$\dashint$ represents the Cauchy principal value (CPV).
For interior points $x_k$ we have $\alpha_k=1$ and CPV integrals are replaced by regular integrals.

In this case
$$
\lambda_{1,k}(\calp)= \left[\int_{\Omega_s^{k}} \frac{1}{\kappa}\rho c p_1 v\, d\Omega, \int_{\Omega_s^{k}} \frac{1}{\kappa}\rho c p_2 v \, d\Omega, \ldots , \int_{\Omega_s^{k}} \frac{1}{\kappa} \rho c p_Q v \, d\Omega \right]^T,
$$
and
$\lambda_{2,k}(\calp)=\alpha_k \lambda^{(1)}_{2,k}(\calp)+\lambda^{(2)}_{2,k}(\calp)+\lambda^{(3)}_{2,k}(\calp)$, where
\begin{align*}
\lambda^{(1)}_{2,k}(\calp)&= \big[p_1(x_k),p_2(x_k),\ldots, p_Q(x_k) \big]^T,   \\
\lambda^{(2)}_{2,k}(\calp)&= -\left[\dashint_{\Gamma_s^{k}} \frac{\partial v}{\partial n}p_1  \, d \Gamma, \dashint_{\Gamma_s^{k}} \frac{\partial v}{\partial n}p_2  \, d \Gamma, \ldots, \dashint_{\Gamma_s^{k}} \frac{\partial v}{\partial n}p_Q  \, d \Gamma   \right]^T, \\
\lambda^{(3)}_{2,k}(\calp)&= \left[\int_{\Omega_s^{k}}\frac{1}{\kappa}\nabla\kappa\cdot \nabla p_1\, v  \, d\Omega, \int_{\Omega_s^{k}}\frac{1}{\kappa}\nabla\kappa\cdot \nabla p_2\, v  \, d\Omega, \ldots, \int_{\Omega_s^{k}}\frac{1}{\kappa}\nabla\kappa\cdot \nabla p_Q\, v  \, d\Omega \right]^T.
\end{align*}
Finally, we have the
time-dependent linear system of equations
\begin{equation}\label{eqfullsys4}
A^{(1)} \frac{\partial}{\partial t} \u(t) +A^{(2)} \u(t) =\b(t),
\end{equation}
where the $k$-th row of $A^{(\ell)}$ is
\begin{align*}
\a_k^{(\mathrm{\ell})}  &= WP (P WP^T)^{-1} \lambda_{\ell,k}(\calp),\quad \ell = 1,2.
\end{align*}
The components of the right-hand side are
$$
b_k(t) = \int_{\Omega_s^{k}}\frac{1}{\kappa(x)}f(x,t)v(x;x_k)  \, d\Omega+\int_{\partial \Omega_s^k\cap\Gamma_N}u_N(x,t) v(x;x_k)\, d\Gamma.
$$
This technique leads to weakly singular integrals which must
be evaluated by special numerical quadratures.

\subsection{DMLPG3/6}
In both MLPG3
and MLPG6,
the trial and test functions come from the same space.
Therefore they are
Galerkin type techniques and should better be called MLG3 and MLG6.
But they
annihilate
the advantages of DMLPG methods
with respect to
numerical integration,
because the integrands include
shape functions.
Thus we ignore DMLPG3/6 in favour of keeping all benefits of DMLPG
methods. Note that MLPG3/6 are also rarely used in comparison to the other MLPG
methods.

\section{Time Stepping}\label{sec-timestep}
To deal with the time variable in meshless methods, some standard
methods
were proposed in the literature.
The
Laplace transform method \cite{sladek-et-al:2003-1,sladek-et-al:2004-1},
conventional finite difference methods such as forward, central and backward
difference schemes are such techniques.
A method which employs the MLS approximation in both time and space domains,
is another different scheme \cite{mirzaei-dehghan:2011-cmes,mirzaei-dehghan:2012-eabe}.

In our case the linear system (\ref{eqalphasys}) turns into
the time--dependent version (\ref{eqfullsys})
coupled with \eqref{dirichlet-sys}
that could, for instance, be solved like any other linear first--order
implicit
Differential Algebraic Equations (DAE)
system. Invoking an ODE solver
on it would be an instance of the Method of Lines.
If a conventional time--difference scheme such as a Crank-Nicolson method
is employed,
if the
time step
$\Delta t$ remains unchanged,
and if $M=N$, then a single
LU decomposition of
the final stiffness matrix
and corresponding backward and forward substitutions
can be calculated once and for all, and then the final solution vector at the
nodes is obtained by a simple matrix--vector iteration.

The classical MLS approximation can be
used as a postprocessing step to obtain the solution at any other point $x\in\Omega$.
\section{Numerical results}\label{sec-num}
Implementation is done using the basis polynomials
\begin{equation*}
\left\{\frac{(x-z)^\beta}{h^{|\beta|}}\right\}_{0\leqslant|\beta|\leqslant m}
\end{equation*}
where $h$
is 
an average mesh-size, and $z$ is a fixed evaluation
point such as a test point or a Gaussian point for integration in
weak--form based techniques. Here $\beta=(\beta_1,\ldots,\beta_d)\in\N_0^d$
is a multi-index and $|\beta|=\beta_1+\ldots+\beta_d$. If
$x=(\schi_1,\ldots,\schi_d)$
then $x^\beta=\schi_1^{\beta_1}\ldots \schi_d^{\beta_d}$. This choice of basis
function,
instead of $\{x^\beta\}_{0\leqslant|\beta|\leqslant m}$, leads to a
well-conditioned
matrix $P WP^T$ in the (G)MLS approximation. The effect of this variation on the
conditioning has been analytically investigated in \cite{mirzaei-et-al:2012-1}.

A test problem is first considered to compare the results of MLPG and DMLPG
methods.
Then a heat conduction problem in functionally graded materials (FGM) for a finite strip
with an exponential spatial variation of material parameters is investigated.
In numerical results, we use the quadratic shifted scaled basis polynomial functions ($m=2$)
in (G)MLS approximation for both MLPG and DMLPG methods. Moreover, the Gaussian weight function
\begin{equation*}
w_j(x)=
\begin{cases}
\frac{\exp(-(\|x-x_j\|_2/c)^2)-\exp(-(\delta/c)^2)}{1-\exp(-(\delta/c)^2)},& 0\leqslant \|x-x_j\|_2\leq \delta, \\
0, &  \|x-x_j\|_2>\delta
\end{cases}
\end{equation*}
where $\delta=\delta_0h$ and $c=c_0 h$ is used.
The parameter 
$\delta_0$ should be large enough to ensure the regularity of
the
moment matrix $PWP^T$ in (G)MLS approximation.
It depends on the degree of
polynomials in use.
Here we put $\delta_0=2m$.
The constant $c_0$ controls the shape of
the weight function and has influence on the stability and
accuracy of (G)MLS approximation. There is no optimal value for this
parameter at hand. Experiments show that $0.4<c_0<1$ lead to more accurate results.

All routines were written using \textsc{Matlab}$^\copyright$
 and run on a Pentium 4 PC with 2.50 GB of Memory and a twin--core 2.00 GHz CPU.

\subsection{Test problem}
Let $\Omega=[0,1]^2\subset\R^2$ and consider
Equations (\ref{govern1})-(\ref{neumanncond}) with $\rho c =2\pi^2$,
$\kappa=1$ and $f(x,t)=0$.
Boundary conditions
using $x=(\schi_1,\schi_2)\in\R^2$ 
are
\begin{align*}
\frac{\partial u}{\partial n} &= 0, \, (\schi_1, \schi_2=0)
\cup (\schi_1, \schi_2=1),\, \schi_1\in[0,1],\\
u &= e^{-t}\cos(\pi \schi_2),\, (\schi_1=0 , \schi_2), \, \schi_2\in[0,1],\\
u &= -e^{-t}\cos(\pi \schi_2),\, (\schi_1=1, \schi_2),\, \schi_2\in[0,1].
\end{align*}
The initial condition is $u(x,0)=\cos(\pi \schi_1)\cos(\pi \schi_2)$, and
$u(x,t)=e^{-t}\cos(\pi \schi_1)\cos(\pi \schi_2)$ is the exact solution.
Let $t_F=1$ and $\Delta t = 0.01$ in
the 
Crank-Nicolson scheme.
A regular node distribution with distance $h$ in both directions is used.
In Table \ref{tb1} the CPU times used by
MLPG1/2/4/5 and DMLPG1/2/4/5 are compared. As we can immediately see, DMLPG methods are
absolutely faster than MLPG methods. There is no significant
difference between MLPG2 and DMLPG2, because they are both collocation techniques and
no numerical integration
is required.
\begin{center}
INSERT TABLE 1
\end{center}

\begin{figure}[hbt]
\begin{center}
\includegraphics[width=12cm]{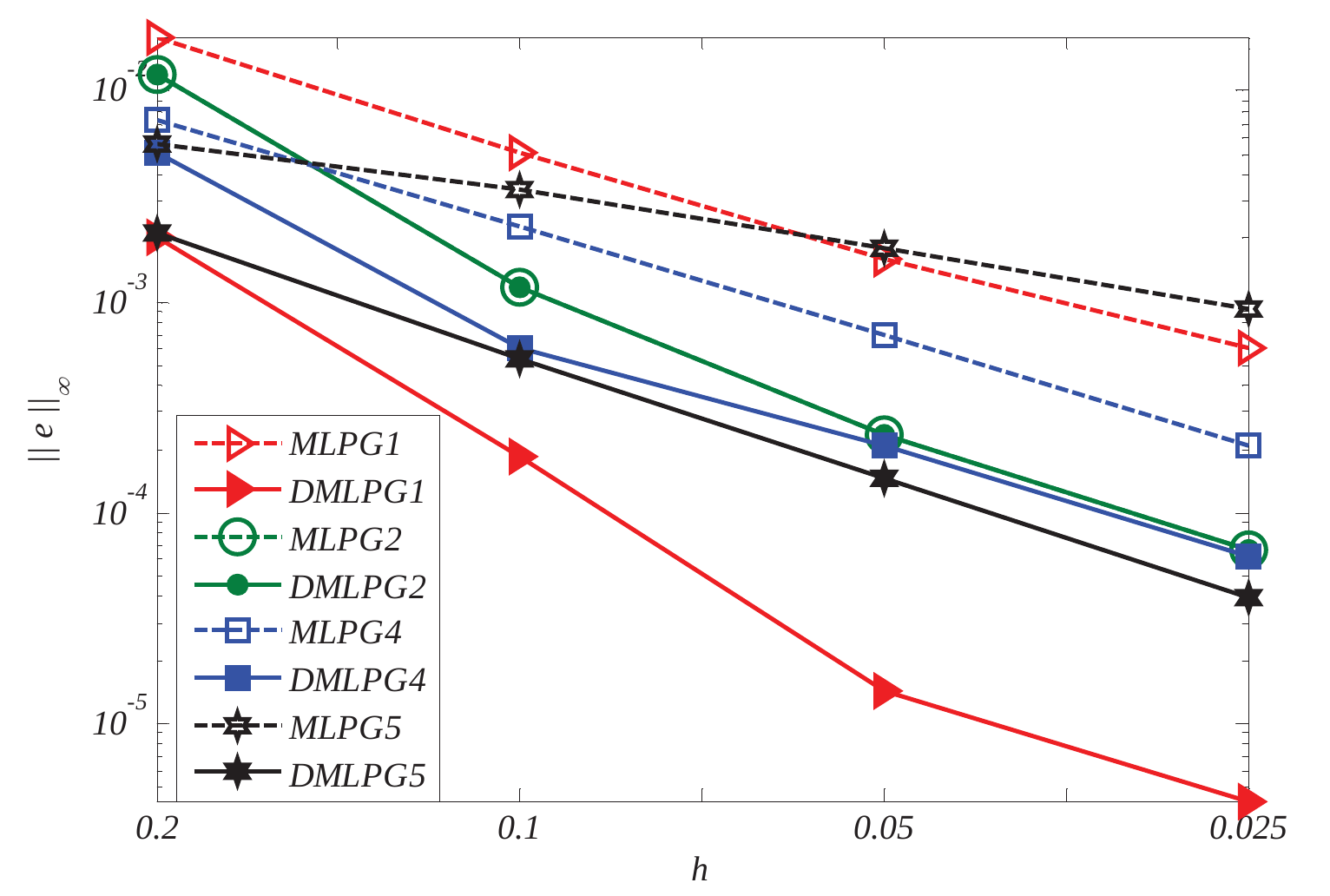}\\
\caption{\small{Comparison of MLPG and DMLPG methods in terms of maximum errors.}}\label{fig1}
\end{center}
\end{figure}

The maximum absolute errors are drawn in
Figure \ref{fig1} and compared. MLPG2 and DMLPG2
coincide, but DMLPG1/4/5 are more accurate than
MLPG1/4/5. DMLPG1 is the most accurate method among all.
Justification needs a rigorous error and stability analysis which is not presented here.
But, according to \cite{mirzaei-schaback:2011-1,mirzaei-et-al:2012-1} and all
numerical results,
we can expect an error behavior like
$\mathcal O(h^{m+1-k})$, where $k$ is the maximal order of derivatives of $u$
involved in the functional, and if numerical integration and time discretization have even smaller errors.

\begin{center}
INSERT FIGURE 1
\end{center}
For more details see the elliptic problems in \cite{mirzaei-schaback:2011-1}
where the ratios of errors of both method types are compared for $m=2,3,4$.
\subsection{A problem in FGMs}
Consider a finite strip with a unidirectional
variation of the thermal conductivity. The exponential spatial
variation is taken
\begin{equation}\label{conductivity}
\kappa(x)=\kappa_0\exp(\gamma \schi_1),
\end{equation}
with $\kappa_0=17~ \mathrm{W m}^{-1~\circ}\mathrm{C}^{-1}$ and $\rho c = 10^6$. This problem
has been considered in \cite{sladek-et-al:2003-1} using
the meshless LBIE method (MLPG4) with Laplace transform in time,
and in \cite{mirzaei-dehghan:2012-eabe,mirzaei-dehghan:2011-cmes}
using MLPG4/5 with MLS approximation for both time and space domains, and
in \cite{wang-et-al:2006-1}
using a RBF based meshless collocation method with time difference approximation.

In numerical calculations,
a square with a side-length $a=4$ cm and a $11\times 11$
regular node distribution is used.

Boundary conditions are imposed as bellow:
the left side is kept to zero temperature and the right side
has the Heaviside step time variation i.e., $u=TH(t)$
with $T=1^\circ C$. On the top and bottom sides the heat flux vanishes.

We employed the
ODE solver \texttt{ode15s} from MATLAB for the final DAE system,
and we used the relative and absolute tolerances \texttt{1e-5} and \texttt{1e-6}, respectively.
With these, we solved on a time interval
of \texttt{[0 60]} with initial condition vector $u_0$ at time $0$.
The Jacobian matrix can be defined in advance because it is constant in our linear DAE.
The integrator will detect stiffness of the system automatically and adjust its local stepsize.

In special case with an
exponential parameter $\gamma=0$ which corresponds to a homogeneous
material the analytical
solution
\begin{align*}
u(x,t)=\frac{T\schi_1}{a}+\frac{2}{\pi}\sum_{n=1}^\infty \frac{T\cos n\pi}{n}\sin\frac{n\pi \schi_1}{a}
\times\exp\left(-\frac{\alpha_0n^2\pi^2t}{a^2}\right),
\end{align*}
 is available. It can be used to check the accuracy of the
present numerical method.

\begin{figure}
\begin{center}
\includegraphics[width=13cm]{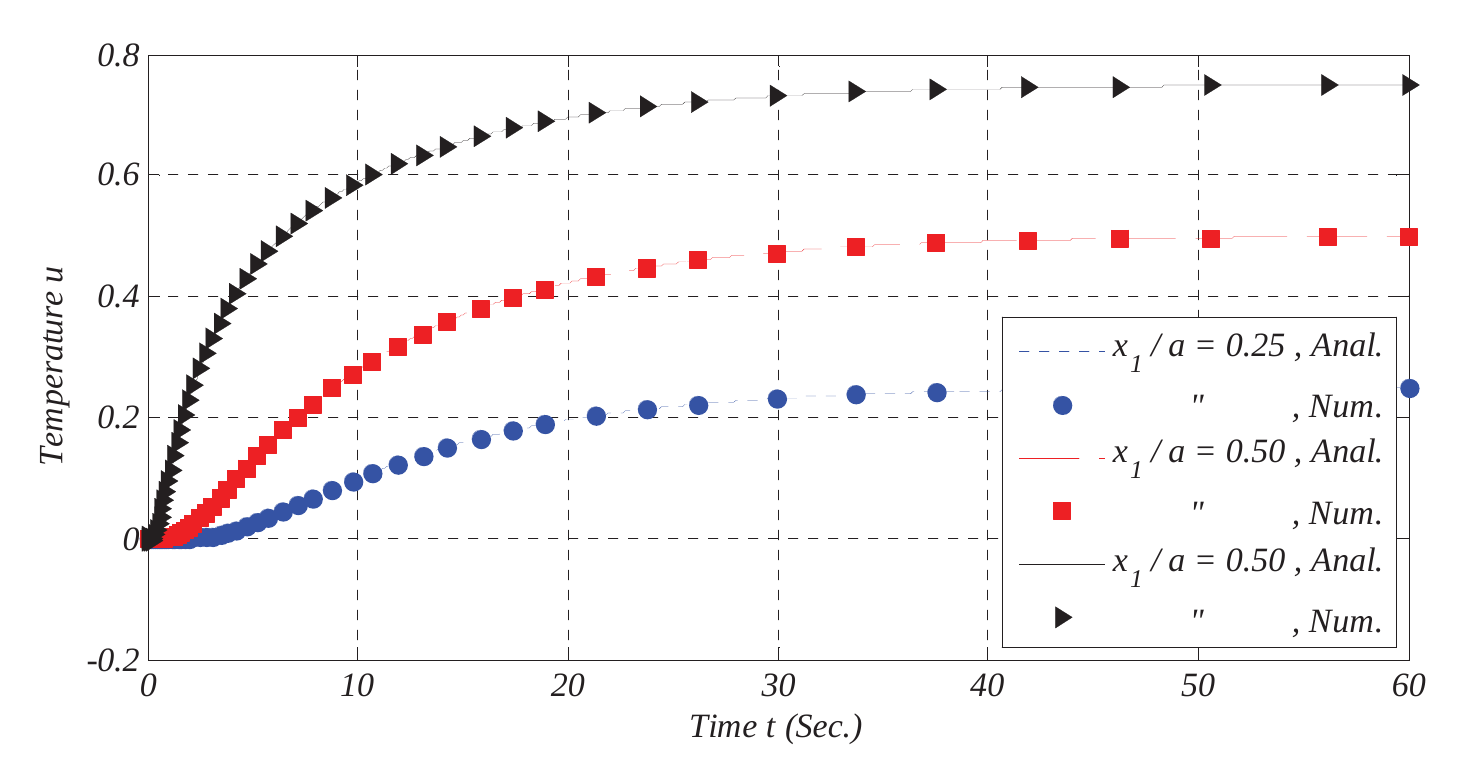}
\caption{\small{Time variation of the temperature at three positions with $\gamma=0$.}}\label{fig2}
\end{center}
\end{figure}

Numerical results
are computed at three locations along the $\schi_1$-axis
with $\schi_1/a=0.25$, $0.5$ and $0.75$.
Results are depicted
in Fig. \ref{fig2}. An excellent agreement between
numerical and analytical solutions is obtained.
\begin{center}
INSERT FIGURE 2
\end{center}

\begin{figure}
\begin{center}
\includegraphics[width=13cm]{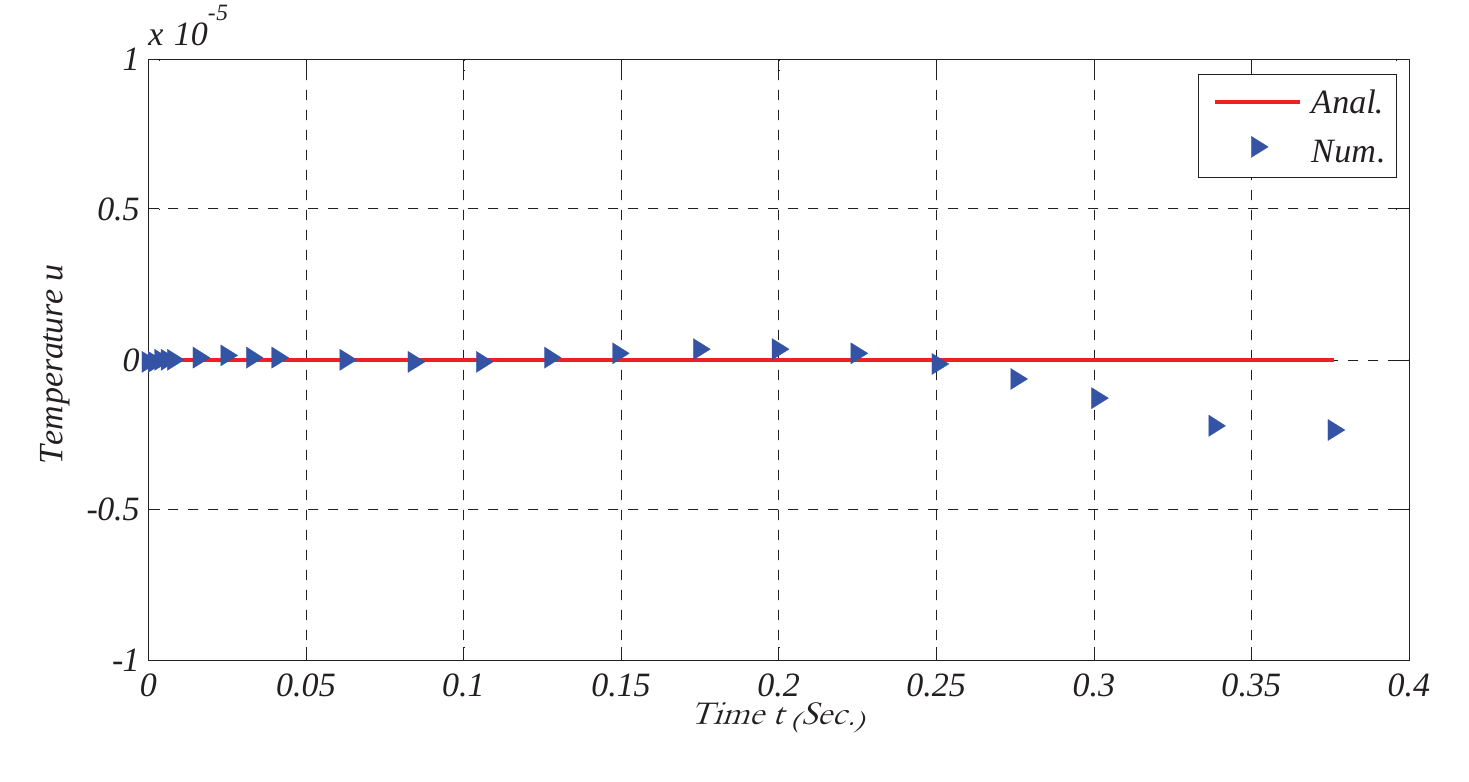}
\caption{\small{Accuracy of method for early time instants at position $x_1/a =0.5$.}}\label{fig3}
\end{center}
\end{figure}

\begin{figure}
\begin{center}
\includegraphics[width=13cm]{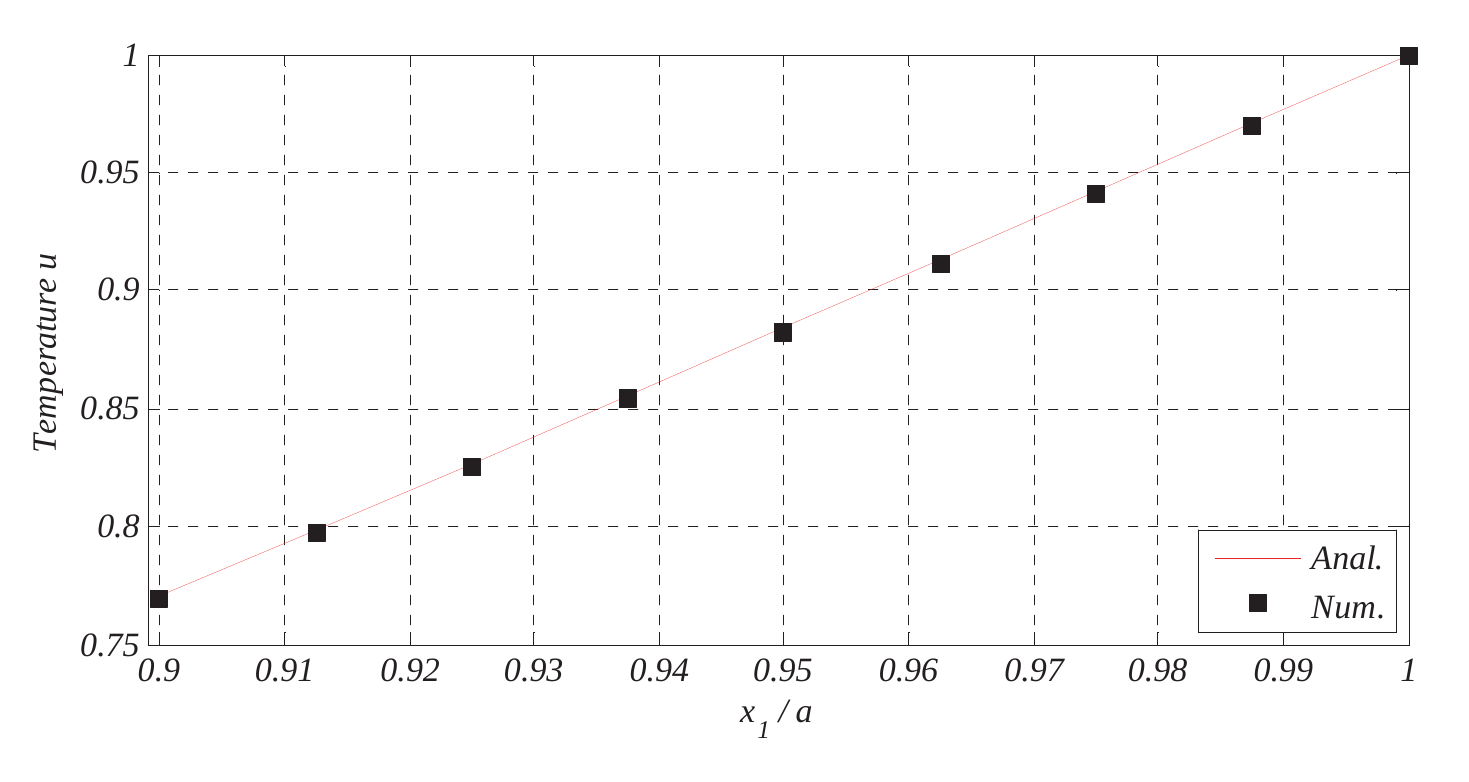}
\caption{\small{Accuracy of method for points close to the thermal shock at time $t(70)$ sec.}}\label{fig4}
\end{center}
\end{figure}

It is known that the numerical results are rather
inaccurate at very early time instants and at
points close to the application of thermal shocks. Therefore in
Fig. \ref{fig3} we have compared the numerical and analytical
solutions at very early time instants ($t\in[0,0.4]$).
Besides, in Fig. \ref{fig4} the numerical and analytical solutions
at points very close to the application of thermal shocks are given and
compared for sample time $t(70)\approx 10.5$ sec.
\begin{center}
INSERT FIGURES 3 AND 4
\end{center}

\begin{figure}
\begin{center}
\includegraphics[width=13cm]{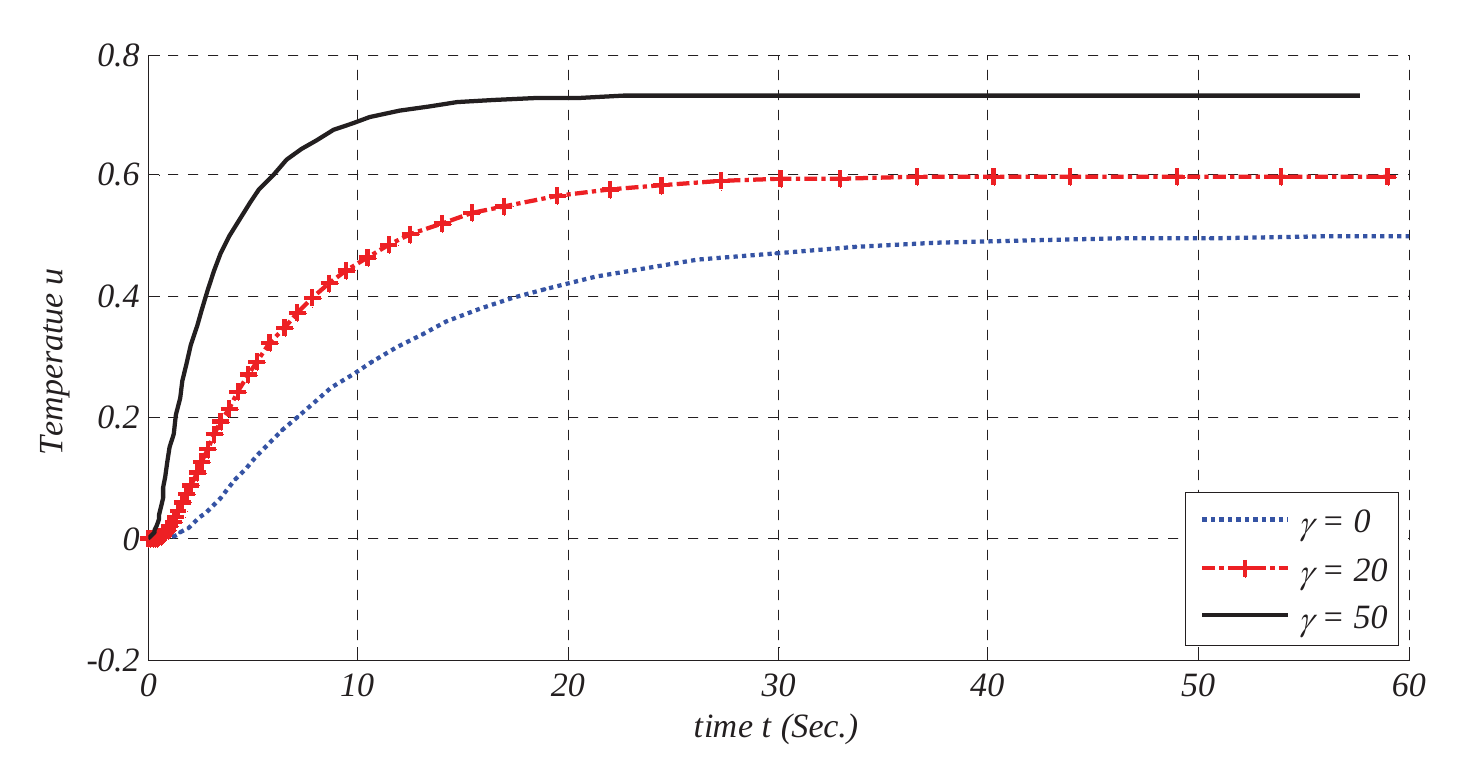}
\caption{\small{Time variation of the temperature at position $x_1/a = 0.5$ for $\gamma = 0,20,50$ m$^{-1}$.}}\label{fig5}
\end{center}
\end{figure}

The discussion above concerns heat conduction in
homogeneous materials
in a case where
analytical solutions
can be used for verification. Consider now
the cases $\gamma = 0$,  $20$, $50$, and $100$ m$^{-1}$, respectively. The
variation of temperature with time for the
three first  $\gamma$-values at position $\schi_1/a = 0.5$  are presented in Fig. \ref{fig5}.
The results are in good agreement with Figure 11
presented in \cite{wang-et-al:2006-1}, Figure 6 presented in
\cite{mirzaei-dehghan:2012-eabe} and Figure 4 presented in \cite{mirzaei-dehghan:2011-cmes}.
\begin{center}
INSERT FIGURE 5
\end{center}

\begin{figure}
\begin{center}
\includegraphics[width=13cm]{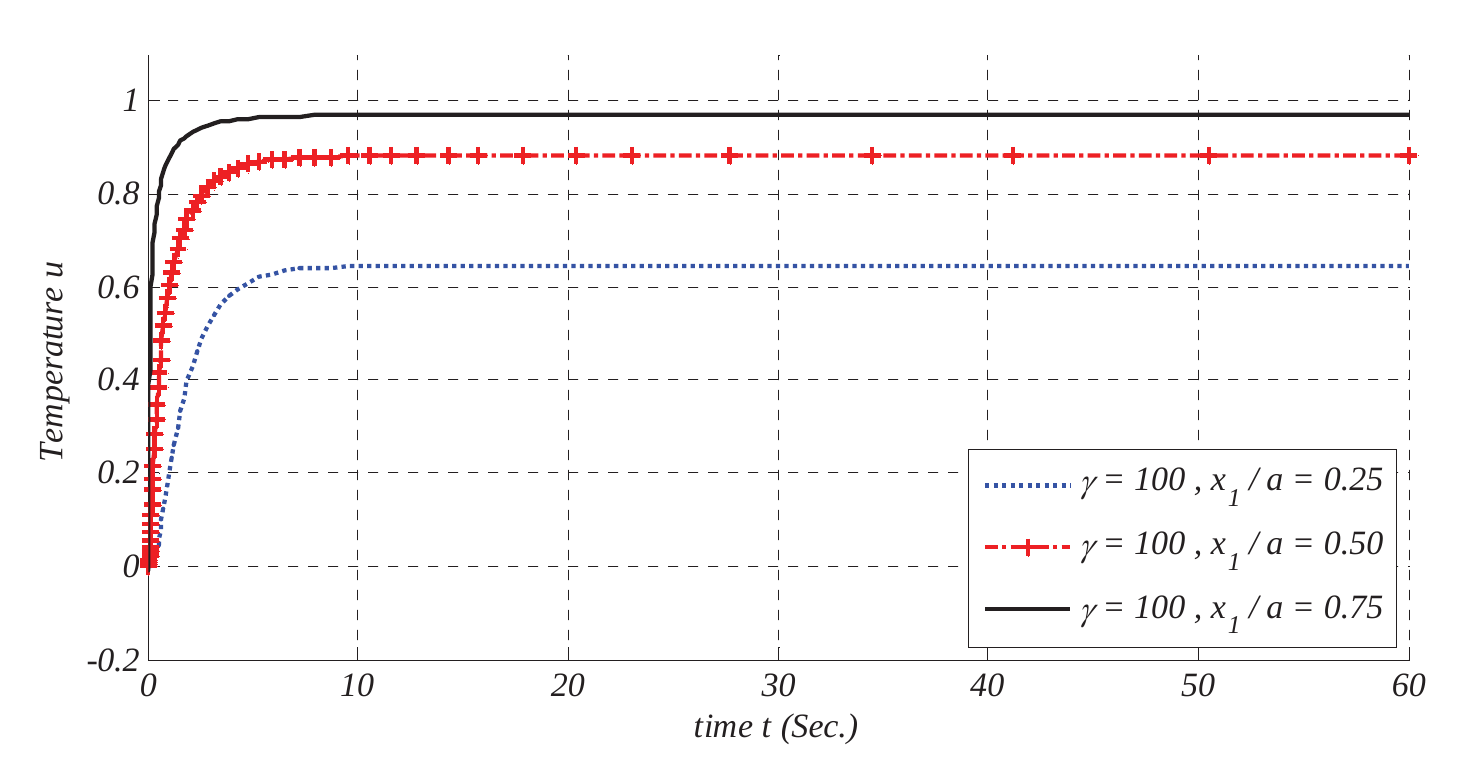}
\caption{\small{Time variation of the temperature at positions $x_1/a =0.25,0.5,0.75$ for $\gamma =100$ m$^{-1}$.}}\label{fig6}
\end{center}
\end{figure}

In addition, in Fig. \ref{fig6} numerical results are
depicted for $\gamma = 100$ m$^{-1}$.
For high values of $\gamma$, the steady state solution is
achieved rapidly.

\begin{center}
INSERT FIGURE 6
\end{center}

It is found from Figs. \ref{fig5} and \ref{fig6} that the temperature
increases with an increase in $\gamma$-values.

\begin{figure}
\begin{center}
\includegraphics[width=13cm]{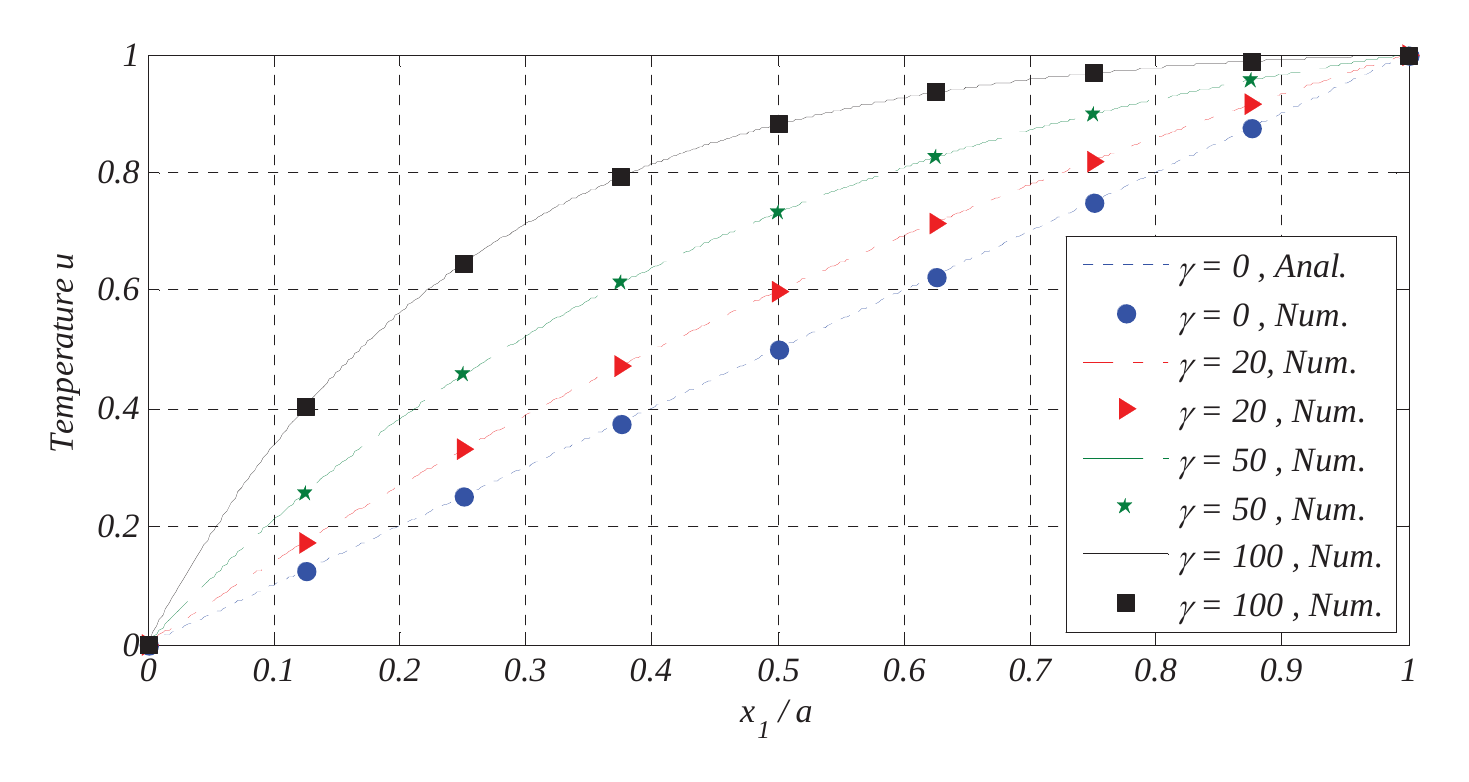}
\caption{\small{Distribution of temperature along $x_1$-axis under steady-state loading conditions.}}\label{fig7}
\end{center}
\end{figure}

For the final steady state, an analytical solution can be
obtained as
\begin{align*}
u(x,t\rightarrow\infty)=T\frac{\exp(-\gamma \schi_1)-1}{\exp(-\gamma a)-1},\quad
 \left(u\rightarrow T\frac{\schi_1}{a}, \mathrm{~as~} \gamma\rightarrow0\right).
\end{align*}
Analytical and numerical results computed at time $t =60$ sec.
are presented in Fig. \ref{fig7}. Numerical results
are in good agreement with analytical solutions for the
steady state temperatures.
\begin{center}
INSERT FIGURE 7
\end{center}

\section*{Acknowledgments}
The first author was financially supported by the Center of Excellence for Mathematics, University of Isfahan.


\bibliographystyle{spmpsci}      


\newpage
\begin{center}
\begin{table}
\begin{tabular}{llllllllllll}
  \hline
               &\fs{Type 1}   &               &&\fs{Type 2}  &             &&\fs{Type 4}   &              &&\fs{Type 5}    &           \\
  \cline{2-3}\cline {5-6} \cline {8-9} \cline {11-12}
  \fs{$h$}     &\fs{MLPG}     & \fs{DMLPG}    && \fs{MLPG}   &\fs{DMLPG}   && \fs{MLPG}    & \fs{DMLPG}   && \fs{MLPG}     & \fs{DMLPG}\\
  \hline

  \fs{$0.2$}   &\fs{$4.3$}    &\fs{$0.2$}     &&\fs{$0.2$}   & \fs{$0.2$}  &&\fs{$1.9$}    & \fs{$0.2$ }  &&\fs{$1.4$}     & \fs{$0.2$}\\
  \fs{$0.1$}   &\fs{$22.6$}   &\fs{$0.3$}     &&\fs{$0.3$}   & \fs{$0.3$}  &&\fs{$9.8$}    & \fs{$0.3$ }  &&\fs{$6.8$}     & \fs{$0.3$}\\
  \fs{$0.05$}  &\fs{$116.4$}  &\fs{$1.4$}     &&\fs{$0.8$}   & \fs{$0.6$}  &&\fs{$52.9$}   & \fs{$1.1$ }  &&\fs{$35.6$}    & \fs{$1.2$}\\
  \fs{$0.025$} &\fs{$855.8$}  &\fs{$9.6$}     &&\fs{$8.3$}   & \fs{$7.0$}  &&\fs{$446.5$}  & \fs{$8.3$ }  &&\fs{$302.2$}   & \fs{$8.5$}\\

  \hline\\
\end{tabular}
\caption{\small{Comparison of MLPG and DMLPG methods in terms of CPU times used (Sec.)}}\label{tb1}
\end{table}
\end{center}



\end{document}